\documentclass[12pt,leqno]{amsart}
\usepackage{fullpage}
\usepackage{amsmath}
\usepackage{amsfonts}
\usepackage{amssymb}
\usepackage{graphicx,xy}

\def\cal{\mathcal}

\def\F2{\mathbb F_2}
\def\A2{{\mathcal A}_2}

\def\N{{\mathbb N}}
\def\Sp{{\cal{S}\rm{p}}}
\def\Op{{\cal{O}\rm{p}}}
\def\TI{{\cal{T}}_I}
\def\TJ{{\cal{T}}_J}
\def\Tor{{\rm Tor}}

\def\Epi{{\rm Epi}}

\def\Bij{{\rm Bij}}

\def\id{{\rm id}}
\def\In{{\rm In}}

\def\Ob{{\rm Ob}}
\def\Sh{{\rm Sh}}

\def\dgvs{{\rm dgvs}}
\def\kfield{{\bf k}}

\newtheorem{thm}{Theorem}[subsection]
\newtheorem{lem}[thm]{Lemma}
\newtheorem{prop}[thm]{Proposition}
\newtheorem{cor}[thm]{Corollary}
\theoremstyle{definition}
\newtheorem{defn}[thm]{Definition}

\newtheorem{nota}[thm]{Notation}
\newtheorem{sub}[thm]{}
\theoremstyle{remark}

\newcommand{\petitcup}{\mathop{\cup}}

\xyoption{all}
\CompileMatrices

\begin{document}

\title[Koszul duality of the category of trees and bar construction for operads]{Koszul duality of the category of trees and bar construction for operads}
\author{Muriel Livernet}
\address{Universit\'e Paris 13, CNRS, UMR 7539 LAGA, 99 avenue
  Jean-Baptiste Cl\'ement, 93430 
  Villetaneuse, France}
\email{livernet@math.univ-paris13.fr}

\date{\today}
\begin{abstract}
In this paper we study a category of trees $\TI$ and prove that it is a Koszul category. Consequences are the interpretation of the reduced bar construction of operads of Ginzburg and Kapranov as the Koszul complex of this category, and the interpretation of operads up to homotopy as a functor from the minimal resolution of $\TI$ to the category of graded vector spaces. We compare also three different bar constructions of operads. Two of them have already been compared by Shnider-Von Osdol and Fresse.
\end{abstract}
\maketitle

\section*{Introduction}

The bar construction is an old machinery that applies to different objects, as algebras,
monads \cite{May72} or  categories \cite{Mitchell72}. Ginzburg and Kapranov built a bar construction $B^{GK}$ for operads in \cite{GinKap94}, as an analogue of the bar construction for algebras. Except that, stricto sensu, this bar construction is not the exact analogue of the bar construction for algebras.  C. Rezk in \cite{RezkPhD}, S. Shnider and D. Von Osdol in \cite{ShnVOs99},  and B. Fresse in \cite{Fr04} considered another bar construction for operads, denoted by $B^{\circ}$, the one viewing an operad $\cal P$ as a monoid in the monoidal (non-symmetric) category of symmetric sequences, with the plethysm as monoidal structure. B. Fresse proved that the associated complex of the two bar constructions are related by an explicit quasi-isomorphism improving the result by S. Shnider and D. Von Osdol who proved that the  two complexes have isomorphic homology.

\medskip

The purpose of this paper is first to give an interpretation of the original bar construction $B^{GK}$ of Ginzburg and Kapranov. In this process we view an operad as a left module over the category of trees $\TI$ (see also \cite{GinKap94}) and build the bar construction $B$  of this category. The crucial point is that this category is Koszul and it is immediate to see that the original bar construction of Ginzburg and Kapranov is precisely the 2-sided Koszul complex of the category with coefficients in the left-module $\cal P$ and the unit right module.
The bar construction of the category of trees is the same as the triple bar construction of the monad arising from the adjunction
$$\text{Operads} \rightleftarrows \text{Symmetric Sequences},$$
where the left adjoint to the forgetful functor is the free operad functor. Symmetric sequences are also known as species-terminology used in the present paper-- or $\mathbb S$-modules.

\medskip

We prove that the inclusion $B^{GK}\rightarrow B$ factors through the quasi-isomorphism described by fresse $B^{GK}\rightarrow B^{\circ}$ where $B^{\circ}$ denotes the bar construction with respects to the monoidal structure. We describe explicitly the quasi-isomorphism $B^{\circ}\rightarrow B$.

Note that  in \cite{ShnVOs99}, S. Shnider and D. Von Osdol interpret the bar construction $B^{\circ}$ as a bar construction of a category, which is not the same as our category of trees. The one used by the authors is the usual category associated to an operad or a PROP (see e.g. \cite{MayTho78}).

\medskip

The advantage of proving that the category $\TI$ is Koszul is that we can provide a smaller resolution of the category $\TI$ than the usual cobar-bar resolution, inspired by the work of B. Fresse in \cite{Frpr09}. We prove that this resolution yields the definition of operads up to homotopy recovering the original definition given by P. Van der Laan in his PhD thesis
\cite{vanderlaanphD}.

\medskip

The plan of the paper is the following. In section 1 we study the category of trees $\TI$, define the two-sided bar construction $B$, the Koszul complex $K$ and prove that the category  is Koszul (theorem \ref{T:TIKoszul}). In section 2 we compare the three bar constructions $B$, $B^{GK}$ and $B^\circ$ and prove the main theorem \ref{S:main}: the factorization of $B^{GK}\rightarrow B$ through the levelization morphism of B. Fresse. Section 3 is devoted to operads up to homotopy.

\medskip

\noindent {\bf Notation.} Let $\kfield$ be a field of any characteristic. The category of differential graded $\kfield$-vector spaces is denoted by $\dgvs$. An object in this category is often called a {\sl complex}.
The symmetric group acting on $n$ elements is denoted by $S_n$. Let $E,F$ be subsets of a set $G$. The notation $E\sqcup F=G$ means that $\{E,F\}$ forms a partition of $G$, that is, $E\cup F=G$ and $E\cap F=\emptyset$. To any set $E$, one associates the vector space $\kfield[E]$  spanned  by $E$.

\medskip

\noindent{\bf Acknowledgments.}  It is a pleasure here to thank the Nankai university and the Chern Institute and the organizers of the summer school and the conference "Operads and universal algebra"
in july 2010. There I had valuable discussion with P.-L. Curien, P. Malbos and Y. Guiraud. This work is inspired by the paper \cite{Frpr09} where B. Fresse. I'd like to thank him for the discussion we had on this subject.


\section{The tree category is Koszul}\label{SS:cat}

\subsection{The tree category $\TI$}

\begin{defn}\label{D:trees} A {\sl tree} is a non-empty connected oriented graph  $t$ with no loops with the property that at each vertex $v$ there is at least one incoming edge and exactly one outgoing edge.
The target of the incoming edges at $v$ is $v$ and the source of the outgoing edge of $v$ is $v$. We allow some edges to have no sources and these edges are called the leaves of the tree $t$. The other ones are called the internal edges of $t$ and we
denote by $E_t$ the set of internal edges of $t$. We denote by $V_t$ the set of vertices of $t$ and
by $\In(v)$ the set of incoming edges (leaves or internal edges) at the vertex $v$. A tree is {\it reduced}
if for every $v\in V_t$ one has $|\In(v)|>1$.

Let $I$ be a set.
An $I$-tree is a tree such that there is a bijection between  the set of its leaves and $I$.
The objects of the category $\TI$ are the reduced $I$-trees. Note that the set of objects of $\TI$ is finite.
Let $t$ be a tree in $\TI$ and $E$ be a set of internal edges which can be empty. 
The tree $t/E$ is the tree obtained by contracting the edges $e\in E$. For a given pair of trees $(t,s)$ the set of 
morphisms  $\TI(t,s)$ is a point if there is $E\subset E_t$ such that $s=t/E$ and is empty if not. Note that if there is
$E\subset E_t$ such that $s=t/E$ then $E$ is unique. The category $\kfield\TI$ is the $\kfield$-linear category spanned by $\TI$:
it has the same set of objects and has for morphisms $\kfield\TI(t,s)=\kfield[\TI(t,s)]$. When $s=t/E$ we denote again by $E$ the basis of the one dimensional vector space
$\kfield\TI(t,t/E)$.

A {\sl left $\TI$-module} is a covariant functor $\TI\rightarrow \dgvs$
and a {\sl right $\TI$-module} is a contravariant functor $\TI\rightarrow \dgvs$. 
To any  left $\TI$-module $L$ and right $\TI$-module $R$ we associate the differential graded vector 
space
$$R\otimes_{\TI} L=\bigoplus_{t\in\TI} R(t)\otimes L(t)/\sim$$
with $f^*(x)\otimes y\sim x\otimes f_*(y)$
whenever $f\in\TI(t,s), x\in R(s)$ and $y\in L(t)$.

Recall that the Yoneda lemma implies the functorial equivalences
$$\kfield\TI(-,s)\otimes_{\TI} L\cong L(s)\quad \text{ and } \quad  R\otimes_{\TI} \kfield\TI(t,-)\cong R(t).$$

\end{defn}

\subsection{Bar construction for the category $\TI$}\label{S:bar_cat}

\begin{defn}\label{D:barTI}\cite{Mitchell72} The {\sl bar construction} (or standard complex in the terminology of Mitchell) of the $\kfield$-linear category $\kfield\TI$ is a simplicial bifunctor
 $\TI^{op}\times \TI\rightarrow \dgvs$ defined by
\begin{multline*}
B_n(\TI,\TI,\TI)(t,s)=\\
\bigoplus_{s_0,\ldots,s_n\in \TI} \kfield\TI(s_0,s)\otimes \kfield\TI(s_1,s_0)\otimes\ldots\otimes \kfield\TI(s_n,s_{n-1})\otimes \kfield\TI(t,s_n)
\end{multline*}

with the simplicial structure given by

$$\begin{array}{llccc}
\textrm{ for }  0\leq i\leq n &  d_i:& B_n(\TI,\TI,\TI)(t,s) &\rightarrow  &B_{n-1}(\TI,\TI,\TI)(t,s)\\
&&a_0\otimes\ldots\otimes a_{n+1}&\mapsto&a_0\otimes\ldots\otimes a_ia_{i+1}\otimes\ldots a_{n+1}, \\
&&&& \\
\textrm{ for }  0\leq j\leq n &  s_j:& B_n(\TI,\TI,\TI)(t,s) &\rightarrow  &B_{n+1}(\TI,\TI,\TI)(t,s)\\
&&a_0\otimes\ldots\otimes a_{n+1}&\mapsto&a_0\otimes\ldots\otimes a_j\otimes 1 \otimes a_{j+1}\otimes\ldots a_{n+1},
\end{array}$$
\end{defn}

To this simplicial bifunctor is associated the usual complex $(B_n,d=\sum_{i=0}^n(-1)^i d_i)_{n\geq 0}$.
If $s=t/E$ then the complex $B_n$ simplifies as

$$B_n(\TI,\TI,\TI)(t,t/E)=\bigoplus_{E_0\sqcup\ldots\sqcup E_{n+1}=E} \kfield(E_0,\ldots,E_{n+1})$$ with

$$d(E_0,\ldots,E_{n+1})=\sum_{i=0}^n (-1)^i (E_0,\ldots,E_i\sqcup E_{i+1},\ldots, E_{n+1}).$$
This complex is augmented by letting $B_{-1}(\TI,\TI,\TI)(t,s)=\kfield\TI(t,s)$.
From B. Mitchell \cite{Mitchell72} this complex is acyclic and
$B_*(\TI,\TI,\TI)\rightarrow \kfield\TI$ is a free resolution of bifunctors.

\begin{sub}{\it Resolution of left and right $\TI$-modules and $\Tor$ functors.}\label{S:res_TI}  

\smallskip
Let $L$ be  a  left  $\TI$-module and $R$ be a right $\TI$-module. The left $\TI$-module
$B(\TI,\TI,\TI)\otimes_{\TI} L$ is denoted by $B(\TI,\TI,L)$, the right $\TI$-module $R\otimes_{\TI} B(\TI,\TI,\TI)$ is denoted by $B(R,\TI,\TI)$ and the differential graded vector space 
$R\otimes_{\TI} B(\TI,\TI,\TI)\otimes_{\TI} L$ is denoted by $B(R,\TI,L)$.

\smallskip

From B. Mitchell one gets that $B(\TI,\TI,L)$ is a free resolution of $L$
in the category of left  $\TI$-modules and $B(R,\TI,\TI)$ is a free resolution of $R$
in the category of right  $\TI$-modules. Consequently

$$H_n B(R,\TI,L)=\Tor_n^{\TI}(R,L).$$
The Yoneda lemma implies that 
the complex computing the $\Tor$ functor has the following form:
$$\begin{array}{clr}
B_n(R,\TI,L)=&\bigoplus\limits_{(t,E\subset E_t)} \bigoplus\limits_{E_1\sqcup\ldots\sqcup E_n=E} 
R(t/E)\otimes (E_1,\ldots,E_n)\otimes L(t), &n\geq 1 \\
&& \\
B_0(R,\TI,L)=&\bigoplus\limits_{t\in\TI} R(t)\otimes L(t).&
\end{array}$$
with the differential
\begin{multline*}
d(x\otimes(E_1,\ldots,E_n)\otimes y)=(E_1)^*(x)\otimes (E_2,\ldots,E_n)\otimes y \\
+\sum_{i=1}^{n-1} (-1)^i x\otimes (E_1,\ldots, E_i\sqcup E_{i+1}\,\ldots,E_n)\otimes y+(-1)^n x\otimes(E_1,\ldots,E_{n-1})\otimes (E_n)_*(y),
\end{multline*}
where $(E_1)^*=R(E_1: t/(E\setminus E_1)\rightarrow t/E)$ and $(E_n)_*=
L(E_n:t\rightarrow t/E_n).$

\end{sub}

\begin{sub}{\it Normalized bar complex.}
Because $B_*$ is  a simplicial bifunctor, one can mod out by the degeneracies
to get the normalized bar complex of the category
$$N_n(\TI,\TI,\TI)(t,t/E)=\bigoplus_{E_0\sqcup\ldots\sqcup E_{n+1}=E,\atop{ E_i\not=\emptyset \text{ for } 1\leq i\leq n}}
 \kfield(E_0,\ldots,E_{n+1})$$
 as well as the normalized bar complexes with coefficients $N_*(R,\TI,\TI), N_*(\TI,\TI,L)$ and $N_*(R,\TI,L)$. 
Furthermore for any left $\TI$-module $L$ and right $\TI$-module $R$ one has  quasi-isomorphisms $$B_*(\TI,\TI,L)\rightarrow N_*(\TI,\TI,L)\rightarrow L$$ in
 the category of left $\TI$-modules and  quasi-isomorphisms
$$B_*(R,\TI,\TI)\rightarrow N_*(R,\TI,\TI)\rightarrow R$$ in
 the category of right $\TI$-modules and quasi-isomorphisms in $\dgvs$
$$B(R,\TI,L)\rightarrow N(R,\TI,L).$$ 

Since $N_*(R,\TI,\TI)$ is a free right $\TI$-module and $N_*(\TI,\TI,L)$ is a free left $\TI$-module
one can either use the bar complex or the normalized bar complex in the sequel, as free 
resolutions of $L$ or $R$.

 \end{sub}

\subsection{The Koszul complex of the category $\TI$}

\begin{nota}\label{N:bt}
 For any tree $t$ we denote by $b_t$ the left and right $\TI$-module which sends $t$ to $\kfield$ and $s\not=t$ to $0$.
If $t$ is the corolla $c_I$ we use the notation $b_I$ instead of $b_{c_I}$.
\end{nota}

Let $E=\{e_1,\ldots,e_n\}$ be a finite set with $n$ elements. 
Let $\kfield[E]$ be the $n$-dimensional vector space spanned by $E$. 
The vector space $\Lambda^n(\kfield[E])$ is a one dimensional vector space. Let  $e_1\wedge\ldots\wedge e_n$ be a basis.

We would like to define the Koszul complex of the category $\TI$ as a bifunctor 
$$K(\TI,\TI,\TI) : \TI^{op}\times \TI\rightarrow \dgvs.$$  
For any pair of trees $(t,s)$, if  there is no $E$ such that $s=t/E$ we let
$K(\TI,\TI,\TI)(t,s)=0$. If $s=t/E$ we let
$$K(\TI,\TI,\TI)(t,t/E)=\\
\bigoplus_{F\sqcup G\subset E} \kfield\TI(t/(F\sqcup G),t/E)\otimes\Lambda^{|G|}(\kfield[G])\otimes \kfield\TI(t,t/F)$$
For any $F \sqcup (G=\{e_1,\ldots,e_g\})\sqcup H=E$ we define
\begin{multline*}
d(H\otimes e_1\wedge\ldots\wedge e_g\otimes F)=
\sum_{i=1}^g (-1)^{i-1} H\cup\{e_i\}\otimes e_1\wedge\ldots\wedge\hat{e_i}\wedge\ldots\wedge e_g\otimes F+ \\
\sum_{i=1}^g (-1)^i H\otimes e_1\wedge\ldots\wedge\hat{e_i}\wedge\ldots\wedge e_g\otimes 
F\cup\{e_i\}.
\end{multline*}

\begin{lem}
The map $d$ satisfies $d^2=0$.
\end{lem}

\begin{proof}The map $d$ splits into two parts $d_l+d_r$. 

One has $d_ld_r+d_rd_l=0$:  if $x_i$ denotes the element $ e_1\wedge\ldots\wedge\hat{e_i}\wedge\ldots\wedge e_g$ and $x_{i,j}$, $i<j$ denotes the element  
$e_1\wedge\ldots\wedge\hat{e_i}\wedge\ldots\wedge\hat{e_j}\wedge\ldots\wedge  e_g$ then

\begin{multline*}
(d_ld_r+d_rd_l)(H\otimes e_1\wedge\ldots\wedge e_g\otimes F)=d_l(\sum_{i=1}^g (-1)^i H\otimes x_i\otimes F\cup\{e_i\})+d_r(\sum_{j=1}^g (-1)^{j-1} H\cup\{e_j\}\otimes x_j\otimes F)\\
=\sum_{j<i}(-1)^{i+j-1}   H\cup\{e_j\}\otimes x_{j,i}\otimes 
F\cup\{e_i\}+\sum_{j>i}(-1)^{i+j}   H\cup\{e_j\}\otimes x_{i,j}\otimes F\cup\{e_i\}+\\
\sum_{i<j}(-1)^{i+j-1} H\cup\{e_j\}\otimes x_{i,j}\otimes F\cup\{e_i\}+\sum_{i>j} (-1)^{i+j-2} H\cup\{e_j\}\otimes x_{j,i}\otimes F\cup\{e_i\}=0.
\end{multline*}
Let $V$ be an n-dimentional vector space. Let $\cal V=\{v_1,\ldots,v_n\}$ be a basis of $V$ with a given order  $v_1<\ldots<v_n$. 
Recall that the Koszul complex
$\Lambda(V)\otimes S(V)$ has the following differential
$$d(x_1\wedge\ldots\wedge x_p\otimes y_1\ldots y_q)=\sum_{i=1}^p (-1)^i x_1\wedge\ldots\wedge \hat{x_i}\wedge\ldots\wedge x_p\otimes x_iy_1\ldots y_q.$$
This complex splits into subcomplexes
$$(\Lambda(V)\otimes S(V),d)=\bigoplus\limits_{\emptyset\not=W\subset \cal V} (C_*^W,d_W)$$
where 
$$C_p^W(V)=\bigoplus_{\{x_1<\ldots<x_p;y_1\leq \ldots\leq y_q\}=W} k[x_1\wedge\ldots\wedge x_p\otimes y_1\ldots y_q].$$
For $F\subset E$ we let $V_F$ be the vector space with basis $\cal V_F=\{e_k,e_k\not\in F\}$. The map $d_l$ corresponds to the differential $d_{\cal V_F}$ of $C_*^{\cal V_F}(V_F)$ and $d_l^2=0$. The same is true for $d_r$, with  $\cal V_H$.
\end{proof}

Note that the Koszul complex is augmented by letting $K_{-1}(\TI,\TI,\TI)(t,s)=\kfield\TI(t,s)$.

\begin{sub}{\it The Koszul complex of the category $\TI$ with coefficients}

\smallskip
  Let $L$ be  a  left  $\TI$-module and $R$ be a right $\TI$-module. The left $\TI$-module
$K(\TI,\TI,\TI)\otimes_{\TI} L$ is a free left $\TI$-module denoted by $K(\TI,\TI,L)$. The right $\TI$-module $R\otimes_{\TI} K(\TI,\TI,\TI)$ is a free right $\TI$-module denoted by 
$K(R,\TI,\TI)$. The differential graded vector space 
$R\otimes_{\TI} K(\TI,\TI,\TI)\otimes_{\TI} L$ is denoted by $K(R,\TI,L)$.

\smallskip
Let $t$ be a tree in $\TI$ and $s=t/E$ for a given $E\subset E_t$. The right $\TI$-module
$K(b_s,\TI,\TI)$ has the following form
$$K(b_s,\TI,\TI)(t)=\bigoplus_{G\sqcup F=E} \Lambda^{|G|}(\kfield[G])\otimes \kfield[F].$$
From the previous proof one gets that it corresponds to a summand of the Koszul complex
$\Lambda(\kfield[E])\otimes S(\kfield[E])$. If $E$ is non empty, this complex is acyclic (see e.g. \cite{Serre65}) and if $E$ is empty it is $\kfield$ in degree $0$. As a consequence we have the theorem
\end{sub}

\begin{thm}\label{T:Kresol} The augmentation $\epsilon:K(b_s,\TI,\TI)\rightarrow b_s$ is a quasi-isomorphism, thus $K(b_s,\TI,\TI)$ is a free resolution of $b_s$ in the  category
of right $\TI$-modules. 

The augmentation $\epsilon: K(\TI,\TI,b_t)\rightarrow b_t$ is a quasi-isomorphism, thus $K(\TI,\TI,b_t)$ is a free resolution of $b_t$ in the  category
of  left $\TI$-modules.
\end{thm}

\subsection{The category $\TI$ is Koszul}\label{S:sub_Koszul}

The aim of this section is to prove that the homology of the complex $N(b_s,\TI,b_t)$ is concentrated in top degree with value $K(b_s,\TI,b_t)$ which amounts to say that the category $\TI$ is Koszul. 

\begin{lem} The map $\kappa:K(\TI,\TI,\TI)\rightarrow B(\TI,\TI,\TI)$ defined by
$$\kappa(t,t/E)(H\otimes e_1\wedge\ldots\wedge e_n\otimes F)=\sum_{\sigma\in S_n} \epsilon(\sigma) H\otimes (e_{\sigma(1)},\ldots,e_{\sigma(n)})\otimes F$$
is a natural transformation of bifunctors.
For any right $\TI$-module $R$  the induced map
$R\otimes_{\TI} \kappa$ commutes with  the augmentation maps $K(R,\TI,\TI)\rightarrow R$ and $B(R,\TI,\TI)\rightarrow R$.

For any left $\TI$-module $L$ the induced map
$\kappa\otimes_{\TI} L$ commutes with the augmentation maps.
\end{lem}

\begin{proof}The only thing we need to prove is that $\kappa$ commutes with the differentials. One has

\begin{align*}
d\kappa(H\otimes e_1\wedge\ldots\wedge e_n\otimes F)=
\sum_{\sigma\in S_n} \epsilon(\sigma)
( H\cup e_{\sigma(1)})\otimes (e_{\sigma(2)},\ldots,e_{\sigma(n)})\otimes F+\\
 \sum_{i=1}^{n-1}(-1)^i\sum_{\sigma\in S_n}\epsilon(\sigma) H\otimes (e_{\sigma(1)},\ldots,e_{\sigma(i)}\cup_{\sigma(i+1)},\ldots,e_{\sigma(n)})\otimes F+\\
(-1)^n \sum_{\sigma\in S_n}\epsilon(\sigma) H\otimes (e_{\sigma(1)},\ldots,e_{\sigma(n-1)})\otimes F\cup e_{\sigma(n)}.
\end{align*}
The middle term vanishes. 
For the first term, we split the sum over $S_n$ into sums over $\sigma\in S_n$
such that $\sigma(1)=i$. Such a $\sigma$ is a composite $\tau\rho$ with $\tau$ having $i$ as fixed point and with $\rho$ being  the cycle $1\rightarrow i\rightarrow i-1\rightarrow\ldots\rightarrow 2\rightarrow 1$.
Hence $\epsilon(\sigma)=\epsilon(\tau)(-1)^{i-1}$. Thus the first term writes
$$\sum_{i=1}^n (-1)^{i-1} \sum_{\tau\in S_n,\tau(i)=i} \epsilon(\tau) H\cup e_i\otimes (e_{\tau(1)},\ldots,\hat e_i,\ldots e_{\tau(n)})\otimes F.$$
For the last term, we split the sum over $S_n$ into sums over $\sigma\in S_n$
such that $\sigma(n)=i$. Such a $\sigma$ is a composite $\tau\eta$ with $\tau$ having $i$ as fixed point and with $\eta$ being the cycle $i\rightarrow i+1\rightarrow \ldots\rightarrow n\rightarrow i$.
Hence $\epsilon(\sigma)=\epsilon(\tau)(-1)^{n-i}$. Thus the last term writes
$$\sum_{i=1}^n (-1)^{i} \sum_{\tau\in S_n,\tau(i)=i} \epsilon(\tau) H\otimes (e_{\tau(1)},\ldots,\hat e_i,\ldots e_{\tau(n)})\otimes F\cup e_i.$$
As a consequence $d\kappa=\kappa d$.
\end{proof}

\begin{prop}The morphisms of right $\TI$-modules
$$\xymatrix{K(b_s,\TI,\TI) \ar[dr]_{\epsilon}\ar[rr]^{b_s\otimes_{\TI}\kappa} &&B(b_s,\TI,\TI)\ar[dl]^{\epsilon} \\
&b_s &}$$
are quasi-isomorphisms.
\end{prop}

\begin{proof}This is a direct corollary of theorem \ref{T:Kresol}
\end{proof}

\begin{thm}\label{T:TIKoszul} The category $\TI$ is Koszul.
\end{thm}

\begin{proof} Because we have quasi-isomorphisms of free left modules
$$K(b_s,\TI,\TI)\rightarrow B(b_s,\TI,\TI)\rightarrow N(b_s,\TI,\TI),$$ 
we have quasi-isomorphisms of differential graded vector spaces
$$K(b_s,\TI,b_t)\rightarrow B(b_s,\TI,b_t)\rightarrow N(b_s,\TI,b_t).$$ 
If $s=t/E$ with $E=\{e_1,\ldots,e_n\}$  then $N(b_s,\TI,b_t)$ is bounded
with top degree $n$. Namely
$$N_n(b_s,\TI,b_t)=\bigoplus\limits_{\sigma\in S_n}  k[(e_{\sigma(1)},\ldots,e_{\sigma(n)})],$$
whereas $K(b_s,\TI,b_t)$ is concentrated in degree $n$, of dimension one with basis $e_1\wedge\ldots\wedge e_n$. Since $K(b_s,\TI,b_t)\rightarrow
N(b_s,\TI,b_t)$ is a quasi-isomorphism, one gets that $K$ is precisely the homology of $N$.
 \end{proof}

 \begin{cor}\label{C:kappa_qiso} For any right $\TI$-module $R$ and left $\TI$-module $L$ the morphisms
 $$\begin{array}{lcccr}
 \kappa:& K(R,\TI,L)&\rightarrow &B(R,\TI,L)& \text{ and } \\
 &&&& \\
 \overline\kappa: & K(R,\TI,L)&\rightarrow& N(R,\TI,L)& 
 \end{array}$$ 
 are quasi-isomorphisms and
 
 $$\Tor^{\TI}_*(R,L)=H_*(K(R,\TI,L)).$$
 \end{cor}

 \begin{proof} It is enough to prove that $\overline\kappa$ is a quasi-isomorphism. Let us consider the filtration by the number of internal vertices
 $$\begin{array}{cc}
 F_p(N)=&\bigoplus\limits_{E:t\rightarrow s, |E|\leq p\atop {E=E_1\sqcup\cdots\sqcup E_n}} R(s)\otimes (E_1,\ldots,E_n)\otimes L(t), \\
 F_p(K)=&\bigoplus\limits_{n\leq p} K_n(R,\TI,L),
 \end{array}$$
 which are subcomplexes of $N(R,\TI,L)$ and $K(R,\TI,L)$ respectively. One has, as complexes,
 $$ \begin{array}{cc}
 F_p(N)/F_{p-1}(N)=&\bigoplus\limits_{E:t\rightarrow s, |E|=p, n\leq p} R(s)\otimes N_n(b_s,\TI,b_t)\otimes L(t), \\
 F_p(K)/F_{p-1}(K)=&\bigoplus\limits_{E:t\rightarrow s, |E|=p} R(s)\otimes K_p(b_s,\TI,b_t)\otimes L(t).
 \end{array}$$
 From theorem \ref{T:TIKoszul} the map $ F_p(K)/F_{p-1}(K)\rightarrow  F_p(N)/F_{p-1}(N)$ is a quasi-isomorphism. Since $F_0(K)=\oplus_s R(s)\otimes L(s)=F_0(N)$, then for every $p$,
 the map $F_p(K)\rightarrow F_p(N)$ is a quasi-isomorphism. This yields the result.
 \end{proof}

 \section{Comparison of three different types of bar constructions for an operad}

 The aim of this section is to compare different kinds of bar constructions for operads, depending on the way we consider operads, either as  left $\TI$-modules, or algebras over the free operad monad, or
 monoids in the monoidal category of species.
 
 Section \ref{S:principle} is an attempt to generalize the bar construction in a framework that applies to  all the cases considered in the paper. Section \ref{S:left_TI} shows that an operad $\cal P$ can be considered as a left  $\TI$-module, yielding to the bar construction $B(R,\TI,\cal P)$ for $R=\TI$ or $R$ a right $\TI$-module. Section \ref{S:cotriple} defines the free operad functor $\cal F$, yielding to the Godemont/May bar construction $B(\cal R,\cal F,\cal P)$ for an $\cal F$-functor $\cal R$. We prove in proposition \ref{P:evalfunctor} that to any right $\TI$-module $R$ is associated an $\cal F$-functor $\pi_I(R)$ such that $B(\pi_I(R),\cal F,\cal P)=B(R,\TI,\cal P)$.  In section \ref{S:bar_mon} we recall the bar construction $B^\circ(R,\cal P,L)$ of an operad $\cal P$ with coefficients in a right $\cal P$-module $R$ and left $\cal P$-module L obtained by viewing an operad as a monoid in
 the monoidal category of species. In section \ref{S:barGK} we recall the original reduced bar construction $B^{GK}$ given by Ginzburg and Kapranov,  which coincides with the Koszul complex
 $K(b_I,\TI,\cal P)$ introduced in section \ref{SS:cat}. We recall the levelization morphism defined by B. Fresse from $B^{GK}$ to $B^{\circ}$. The last section \ref{S:main} is devoted to the factorization
 of $\bar\kappa:K(b_I,\TI,\cal P)\rightarrow N(b_I,\TI,\cal P)$ introduced in section \ref{S:sub_Koszul} through the levelization morphism.

 \subsection{Principle of the bar construction with coefficients}\label{S:principle}\hfill\break

 The paragraph \cite[section 2.3]{MSS02} of the book of M. Markl, S. Shnider and J. Stasheff can serve as our definition of two-sided bar construction. The idea is to work in a ``context'' for which any object $X$ admits the 
notions of left $X$-modules and right $X$-modules, as

\begin{enumerate}
\item A $\kfield$-algebra $X$ with its usual notions of left $X$-module and right $X$-module;
\item a linear category $X$ where left $X$-modules and right $X$-modules are covariant and contravariant functors $X\rightarrow \dgvs$;
\item A monoid $X$ in a monoidal category $(\cal C,\otimes, I)$ where left  $X$-modules $L$ and right 
 $X$-modules $R$ are
objects in $\cal C$ together with maps $X\otimes L\rightarrow X$ and $R\otimes X\rightarrow X$ commuting with the monoid structure of $X$;
\item A monad $X:\cal T\rightarrow \cal T$ where left modules $L$ and right modules $R$
are functors $L:\cal D\rightarrow \cal T$ and $R:\cal T\rightarrow \cal E$ together with natural transformations $\rho:XL\Rightarrow L$ and $\lambda:RX\Rightarrow R$ commuting with the monad structure.
\end{enumerate}

Note that the last example is very close to a monoid in a monoidal category except that $L$ and $R$ are not objects in the same category as $X$.  Certainly the right notion in order to unify all the examples enumerated above is to start with a monoidal category $\cal C$,  a left and right module categories
$\cal L$ and $\cal R$ (see \cite{Ostrik03} for the definition), and pick a monoid $X$ in $\cal C$ and left module $L\in \cal L$ and right module $R\in\cal R$. 

In this context the above examples resume to

\begin{enumerate}
\item The category $\cal C$ is the category of $\kfield$-modules with the tensor product as monoidal structure and $\cal L=\cal C=\cal R$.
\item The category $\cal C$ is the category of bifunctors $X^{op}\times X\rightarrow \dgvs$ with the tensor product defined in section \ref{SS:cat}. The category $\cal L$ is the category of covariant functors $X\rightarrow\dgvs$ 
and $\cal R$ the one of contravariant functors.
\item The category $\cal C$ is the monoidal category $(\cal C,\otimes, I)$ and $\cal L=\cal C=\cal R$
\item The category $\cal C$ is the category of endo-functors $\cal T\rightarrow \cal T$, with the composition as monoidal structure and the category $\cal L$ is the category of functors $\cal D\rightarrow \cal T$ and the category $\cal R$ is the category of functors $\cal T\rightarrow \cal E$. 
\end{enumerate}
 
 \begin{defn}
 We say that a simplicial complex $B_*(R,\cal P,L)$ endowed with an augmentation $\epsilon:B(R,\cal P,L)\rightarrow B_{-1}(R,\cal P,L)$ satisfies the {\sl principle of the simplicial bar construction with coefficients} if
 
 \begin{itemize}
 \item  $\forall n$, $B_n(R,\cal P,\cal P)$ is a free right $\cal P$-module and $\epsilon:B(R,\cal P,\cal P)\rightarrow 
B_{-1}(R,\cal P,\cal P)=R$ is a quasi-isomorphism.
  \item $\forall n$, $B_n(\cal P,\cal P,L)$ is a free left $\cal P$-module and 
$\epsilon:B(\cal P,\cal P,L)\rightarrow B_{-1}(\cal P,\cal P,L)=L$ is a quasi-isomorphism.
 \end{itemize}
 Since we are working in linear categories, the normalized complex $N_*(R,\cal P,L)$ makes sense and we say that it satisfies the {\sl principle of the  bar construction with coefficients} if it satifies the  properties as the ones stated above. More generally a complex $K_*(R,\cal P,L)$ satisfies the {\it principle
 of the bar construction with coefficients} if it satisfies these properties.
 \end{defn} 
 
 The result of section \ref{SS:cat} can be summed up in the following proposition.
 
 \begin{prop} The standard resolution $B_*(R,\TI,L)$ satisfies the principal of the simplicial bar construction. The normalized complex $N_*(R,\TI,L)$ satisfies the principle of the bar construction as well as the Koszul complex $K_*(R,\TI,L)$.
 \end{prop}

 \subsection{Operads as left $\TI$-modules}\label{S:left_TI}\hfill\break

In this section, we recall that operads can be considered as left $\TI$-modules as presented in \cite[section 1.2]{GinKap94}. We will be concerned in the sequel with connected operads ($\cal P(0)=0$ and $\cal P(1)=\kfield$).

 \begin{defn}\label{D:operad} Let $\Bij$ be the category whose objects are finite sets and morphisms are bijections. A {\sl vector species} 
 is a contravariant functor $\cal M: \Bij\rightarrow \dgvs$.  An {\sl operad} is a vector species $\cal P$ together with partial composition maps
 $$\circ_i:\cal P(I)\otimes \cal P(J)\rightarrow \cal P(I\setminus\{i\}\sqcup J), \ \forall i\in I, $$
 and unit $\kfield\rightarrow \cal P(\{x\})$ satisfying functoriality, associativity and unit axioms. 
 A {\sl connected species} is a species $\cal M$ such that $\cal M(\emptyset)=0$ and $\cal M(\{x\})=\kfield$.  We denote by $\bar{\cal M}$ the species 
 $$\begin{cases} \bar{\cal M}(I)=0, & \text{ if } |I|\leq 1 \\
 \bar{\cal M}(I)=\cal M(I), & \text{ if } |I|>1. \end{cases}$$
 A {\sl connected operad} is a connected species which is an operad.
 Let $\Sp$ denote the category of connected species and $\Op$ the category of connected operads.
 \end{defn}

 Let $t$ be a tree in $\TI$ and let  $\cal M$ be a vector species. The graded vector space $\cal M(t)$ is defined by
 
\begin{align}
\cal M(t)=\bigotimes_{v\in V_t} \cal M(\In(v)).
\end{align}

When $\cal P$ is an operad, this definition extends to morphisms in $\TI$ so that one gets
a functor $\cal P:\TI\rightarrow \dgvs$, as follows. Let $e\in E_t$ be an internal edge of $t$ going from $w$ to $v$. By reordering the terms in the tensor product one gets 
$$\cal P(t)= \cal P(\In(v))\otimes \cal P(\In (w))\otimes\underbrace{\otimes_{z\in V_t\setminus\{v,w\}}\cal P(\In(z))}_{X_{v,w}}$$ 
and $\cal P(t\rightarrow t/e):\cal P(t)\rightarrow \cal P(t/e)$ is defined as

$$\circ_e\otimes X_{v,w}: \cal P(\In(v))\otimes \cal P(\In (w))\otimes X_{v,w}  \rightarrow\cal P(\In(v)\setminus\{e\}\sqcup \In (w))\otimes X_{v,w}.$$

Iterating the process, and because of the axioms of the operad,  to any  $E\subset E_t$ 
one gets a well defined map $\cal P(t\rightarrow t/E): \cal P(t)\rightarrow \cal P(t/E)$. Consequently $\cal P$ is a left $\TI$-module. In the sequel we will use the notation $E_*$ for the map $\cal P(t\rightarrow t/E)$.

\medskip

{\it In the sequel we will consider the two-sided bar construction $B(\TI,\TI,\cal P)$ and
$B(R,\TI,\cal P)$ for $\cal P$ an operad considered as a left $\TI$-module and $R$ a right $\TI$-module.}

\subsection{Two-sided bar construction from the free operad functor}\label{S:cotriple} \hfill\break


In \cite[chapter 9]{May72}, May defines $B_\ast(R,C,X)$ for any monad $C$, a $C$-algebra 
$X$ and a $C$-functor $R$ to be $RC^nX$ in degree $n$ with the obvious faces and degeneracies corresponding to the $C$-structure, which satisfies the principle of the simplicial bar construction. The idea generalizes the Godement resolution associated to a triple and constructions used by J. Beck. P. May applied this simplicial resolution to the operad $\cal C_n$ of little $n$-cubes. In \cite{BerMoe06} C. Berger and I. Moerdijk compare this construction for operads with the Boardman-Vogt $W$ construction.

In this section we use this construction and compare it to the bar construction for the category $\TI$, in the spirit of E. Getzler and M. Kapranov in \cite[2.17]{GetKap98}.

Let $C:\cal C\rightarrow \cal  C$ be a monad with structural maps $\mu: C^2\rightarrow C$ and $\eta: \id_\cal C\rightarrow C$. A  {\sl $C$-functor $R$} is a functor $R:\cal C\rightarrow \cal D$ together with a natural transformation $\lambda: RC\Rightarrow R$ satisfying the following identities
\begin{align*}
\lambda\circ R\eta=\id : R\Rightarrow R \\
\lambda\circ R\mu=\lambda\circ \lambda C: RC^2\Rightarrow R
\end{align*}

\begin{defn}
The forgetful functor $\Op\rightarrow\Sp$ admits a left adjoint functor, the {\it free operad functor} 
$$\begin{array}{ccc}
\Sp& \rightarrow & \Op \\
\cal M& \mapsto &\cal F\cal M: I\mapsto\bigoplus\limits_{t\in \TI}\cal M(t)
\end{array}$$ 

The partial composition maps $\circ_i:\cal F(\cal M)(I)\otimes \cal F(\cal M)(J)\rightarrow \cal F(\cal M)(I\setminus\{i\}\sqcup J)$ correspond to the grafting of the root of a tree $s\in \TJ$ on the leave $i$ of a tree $t\in\TI$. When $|I|=1$ we let $\cal F(\cal M)(I)=\kfield$. 

An element in $\cal M(t)\subset \cal F(\cal M)(I)$ writes $(t,E_t,m_t)$.  There is
an injection of species $\cal M\rightarrow \cal F(\cal M)$  where the map $\cal M(I)\rightarrow \cal F(\cal M)(I)$ sends $m$ to 
$(c_I,\emptyset,m)\in \cal M(c_I)$, then identifying  $\cal M(I)$ with $\cal M(c_I)$.

\end{defn}

\begin{sub}{\it The two-sided bar construction.} The free operad functor yields a monad on $\Sp$ denoted also by $\cal F$.  
The tripleability theorem implies that $\cal F$-algebras are exactly operads 
\cite[Theorem 1.2]{GetJon94}.
We denote by $\cal F^{(n)}$ the n-th iteration of $\cal F$. An element in $\cal F^{(n)}(\cal P)(I)$ writes
$(t;E_1,\ldots,E_n,p_t)$ with $t\in\TI$, $E_1\sqcup\ldots\sqcup E_n=E_t$ and $p_t\in\cal P(t)$. The counit $\epsilon$ of the adjunction corresponds to the composition in the left $\TI$-module $\cal P$, namely
$$\begin{array}{lccc}
\epsilon:& \cal F(\cal P) & \rightarrow & \cal P \\
& (t,E_t,p_t) & \mapsto & (E_t)_*(p_t)
\end{array}$$
where $(E_t)_*(p_t)$ is in the component $\cal P(c_I)$ of $\cal F(\cal P)(I)$ that we identify with $\cal P(I)$.
The two-sided bar construction $B_n(\cal F,\cal F, P)(I)$ is the simplicial differential graded vector space $\cal F^{n+1}(\cal P)(I)$ with faces
$d_i: \cal F^{(n+1)}(\cal P)(I)\rightarrow \cal F^{(n)}(\cal P)(I)$
 defined by 
\begin{align*}
d_i(t;E_0,\ldots,E_n,p_t)=&(t;E_0,\ldots,E_i\cup E_{i+1},\ldots,E_n,p_t),\ 0\leq i\leq n-1 \\
d_n(t;E_0,\ldots,E_n,p_t)=&(t/E_n;E_0,\ldots,E_{n-1},(E_n)_*(p_t))
\end{align*}
\end{sub}

As a consequence, comparing with the construction in \ref{S:bar_cat}, we have

\begin{prop} The two-sided bar construction $B(\cal F,\cal F,\cal P)(I)$ coincides with \\
 $B(\TI,\TI,\cal P)(c_I)=B(b_I,\TI,\cal P)$.
\end{prop}

\begin{sub}{\it Right $\TI$-modules and $\cal F$-functors.} Let $R:\TI\rightarrow \dgvs$ be a right $\TI$-module. The functor
$$\begin{array}{cccc}
\pi_I(R):&\Sp&\rightarrow& \dgvs \\
& M &\mapsto &\bigoplus\limits_{t\in\TI} M(t)\otimes R(t)
\end{array}$$
determines an $\cal F$-functor.  In order to define the structural map $\lambda: \pi_I(R)\cal F\Rightarrow \pi_I(R)$ one needs to describe, for any $M\in\Sp$, the map
$\lambda_M:\pi_I(R)\cal F(M)\rightarrow \pi_I(R)(M)$.
The vector space $\pi_I(R)\cal F(M)$ is the direct summand of the vector spaces $\cal F(M)(t)\otimes R(t)$, for $t\in\TI$. An element
in $\cal F(M)(t)$ writes $(t',E_1,E_2,m_{t'})$ with $t'/E_2=t$, $E_1\sqcup E_2=E_{t'}$ and $m_{t'}\in M(t')$. The map
$\lambda_M$ assigns the element $m_{t'}\otimes (E_2)^*(r_t)\in M(t')\otimes R(t')$ 
to the element $(t',E_1,E_2,m_{t'})\otimes r_t\in \cal F(M)(t)\otimes R(t)$.

As an example, the $\cal F$-functor $\pi_I(b_I)$, where $b_I$ has been defined in \ref{N:bt}, 
is the  functor $M\mapsto M(I)$ with structural map 
$$\begin{array}{cccc}
\lambda_M:& \cal F(M)(I)& \rightarrow &M(I) \\
& (t,E_t,m_t) &\mapsto &\begin{cases} 0,& \text{ if } t\not=c_I \Leftrightarrow E_t\not=\emptyset, \\
m_{c_I}, & \text{ if } t=c_I. \end{cases}
\end{array}$$

Comparing with the construction in section \ref{S:bar_cat} one gets easily

\begin{prop}\label{P:evalfunctor}Let $R$ be a right $\TI$-module and let $\cal P$ be an operad. The two-sided bar construction 
$B(\pi_I(R),\cal F,\cal P)$ coincides with $B(R,\TI,\cal P)$.
\end{prop}

 \end{sub}

 \subsection{The bar construction with respect to the monoidal structure $\circ$}\label{S:bar_mon}\hfill\break
 
 As pointed out in the introduction, C. Rezk, S. Shnider, D. Von Osdol and B. Fresse have considered a bar construction for operads related to the fact that operads are monoids in the monoidal category of species, adapting the usual bar construction for algebras. Though the category of species is not monoidal symmetric and the monoidal structure is left distributive with respect to the coproduct but not right distributive, one can still perform the bar construction and then define cohomology theories. For the reader interested by this aspect, we refer to the paper by H.-J. Baues, M. Jibladze and A. Tonks \cite{BJT97}.

 \medskip

The category of connected species admits a monoidal structure given by

$$(\cal M\circ \cal N)(J)=\bigoplus_{J_1\sqcup\ldots\sqcup J_r=J} \cal M(\{1,\ldots,r\})\otimes_{S_r} \cal N(J_1)\otimes\ldots\cal N(J_r),$$
with unit
$$\cal I(J)=\begin{cases} $\kfield$,& \text{ if } |J|=1, \\
 0,& \text{ if } |J|\not=1. \end{cases}$$

A connected operad as defined in
definition \ref{D:operad} is exactly a monoid in the monoidal category of connected species $(\Sp,\circ,\cal I)$. Let $\cal P$ be a connected operad. In the sequel we will use the notation $\cal P(n)$ for $\cal P(\{1,\ldots,n\})$ and
$u(v_1,\ldots,v_k)$ for the image of the element $u\otimes v_1\otimes\ldots\otimes v_k\in \cal P(k)\otimes\cal P(I_1)\otimes\ldots\otimes \cal P(I_k)$ under the structure map $\cal P\circ\cal P\rightarrow \cal P$.

 There exists a simplicial bar construction, a normalized bar construction, and construction with coefficients related to the monoidal structure.
 
 \begin{defn}[\cite{Fr04}] Let $\mathcal P$ be an operad, let $R$ be a right $\mathcal P$-module, that is, a species together with a right action $R\circ \mathcal P\rightarrow R$ satisfying the usual associativity and unit condition of a right module, and let $L$ be a left $\cal P$-module.
 The bar construction with coefficients $R$ and $L$ is the simplicial species
 
 $$B_n^{\circ}(R,\cal P,L)=R\circ\underbrace{\cal P\circ\ldots\circ\cal P}_{n \text{ terms }}\circ L$$
 where faces $d_i$ are induced either by the multiplication $\gamma_{\cal P}:\cal P\circ\cal P\rightarrow \cal P$ or by the left and right action and where degeneracies are induced
 by the unit map $\cal I\rightarrow \cal P$. Modding out by the degeneracies, one gets the normalized bar complex $N^{\circ}(R,\cal P,L)$.
 \end{defn}

 \begin{thm}[\cite{Fr04}] The simplicial complex $B_*^\circ(R,\cal P,L)$ satisfies the principle \ref{S:principle} of the simplicial bar construction with coefficients.
 \end{thm}
 
As pointed out by B. Fresse,  $N_*^{\circ}(R,\cal P ,\cal P)$ is a free resolution of the right $\cal P$-module $R$, but
$N_*^{\circ}(\cal P,\cal P ,L)$  is not. So $N_*^{\circ}(R,\cal P ,\cal P)$ satisfies the "right" principle of the bar construction only.

The species $\cal I$ is a right and left module for any connected operad $\cal P$, using the augmentation map $\epsilon: \cal P\rightarrow \cal I$: 
$$\xymatrix{\cal I\circ\cal P\ar[r]^{\id\circ\epsilon} & \cal I\circ \cal I \ar[r]^{\gamma_{\cal I}}& \cal I & \text{ and }&
\cal P\circ\cal I\ar[r]^{\epsilon\circ\id} & \cal I\circ \cal I \ar[r]^{\gamma_{\cal I}}&\cal I.}$$

In the sequel we will be interested by the bar construction $B_*^\circ(\cal I,\cal P,\cal I)$ with coefficients in the
$\cal P$-module $\cal I$ and its normalized complex $N_*^\circ(\cal I,\cal P,\cal I)$.
An element in $B_{n}^\circ(\cal I,\cal P,\cal I)=\cal P^{\circ n}$ is represented by a tree with $n$ levels as in \cite[section 4.3.1]{Fr04}. As an example the tree 

$$t=\qquad \vcenter{\xymatrix@M=0pt@R=4mm@C=4mm{ & \ar[d] &
\ar[dr] & \ar[d] & \ar[dl] & \ar[dr] & & \ar[dl] \\
*+<3mm>{2}\ar@{.}[r] & *+<3mm>[o][F]{1}\ar@{.}[rr]\ar[dr] & & *+<3mm>[o][F]{w}\ar@{.}[rrr]\ar[dl] & & &
*+<3mm>[o][F]{z}\ar@{.}[r]\ar[d] & \\
*+<3mm>{1}\ar@{.}[rr] & & *+<3mm>[o][F]{v}\ar@{.}[rrrr]\ar[drr] & & & & *+<3mm>[o][F]{1}\ar@{.}[r]\ar[dll] & \\
*+<3mm>{0}\ar@{.}[rrrr] & & & & *+<3mm>[o][F]{u}\ar@{.}[rrr] & & & \\
 }}$$

lives in $\cal P\circ \cal P\circ \cal P=B_3^\circ(\cal I,\cal P,\cal I)$ and has 3 levels. The differential of $t$ is a sum of trees with 2 levels in $\cal P\circ\cal P$.

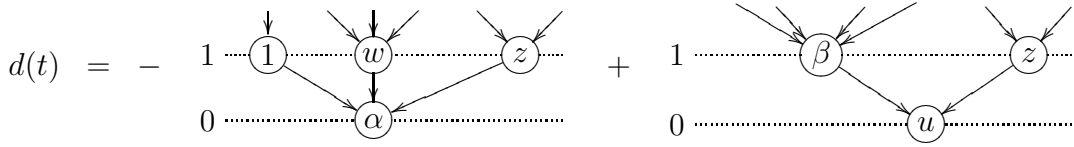
$$\begin{array}{cccccc}
d(t)&=&-&\vcenter{\xymatrix@R=3mm@C=3mm{ & \ar[d] &\ar[dr] & \ar[d] & \ar[dl] & \ar[dr] & & \ar[dl] \\
*+<3mm>{1}\ar@{.}[r] & *+<3mm>[o][F]{1}\ar@{.}[rr]\ar[drr] & & *+<3mm>[o][F]{w}\ar@{.}[rrr]\ar[d] & & &
*+<3mm>[o][F]{z}\ar@{.}[r]\ar[dlll] & \\
*+<3mm>{0}\ar@{.}[rrr] & & &  *+<3mm>[o][F]{\alpha}\ar@{.}[rrrr] & & & & 
 }}&+&
\vcenter{\xymatrix@R=3mm@C=3mm{ & \ar[drr] &
\ar[dr]& & \ar[dl] & \ar[dll] & \ar[dr] & & \ar[dl] \\
*+<3mm>{1}\ar@{.}[rrr] & & & *+<3mm>[o][F]{\beta}\ar@{.}[rrrr]\ar[drr] & & & &  
*+<3mm>[o][F]{z}\ar@{.}[r]\ar[dll] & \\
*+<3mm>{0}\ar@{.}[rrrrr] & & & & & *+<3mm>[o][F]{u}\ar@{.}[rrr] & & & 
 }}
 \end{array}
$$
with $\alpha=u(v,1)\in \cal P(3)$ and $\beta=v(1,w)\in\cal P(4)$.

Note that an element in $N_{n}^\circ(\cal I,\cal P,\cal I)$ is represented by a tree with $n$ levels with the condition that at each level there is at least one vertex labelled by an element in 
$\cal P(r), r\geq 2$. For instance the tree with $3$ levels

$$
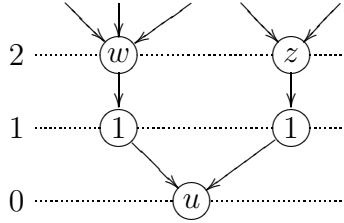
\xymatrix@M=0pt@R=4mm@C=4mm{ & 
\ar[dr] & \ar[d] & \ar[dl] & \ar[dr] & & \ar[dl] \\
*+<3mm>{2}\ar@{.}[rr]& &  *+<3mm>[o][F]{w}\ar@{.}[rrr]\ar[d] & & &
*+<3mm>[o][F]{z}\ar@{.}[r]\ar[d] & \\
*+<3mm>{1}\ar@{.}[rr] & & *+<3mm>[o][F]{1}\ar@{.}[rrr]\ar[dr] & & &  *+<3mm>[o][F]{1}\ar@{.}[r]\ar[dll] & \\
*+<3mm>{0}\ar@{.}[rrr] & & &  *+<3mm>[o][F]{u}\ar@{.}[rrr] & & & \\
 }$$
is zero in $N_{3}^\circ(\cal I,\cal P,\cal I)$.

\subsection{The classical bar construction of operads, and the levelization morphism}\label{S:barGK}
\hfill\break
 
 Ginzburg and Kapranov  introduced in \cite{GinKap94} the reduced bar construction, based on partial compositions, as defined in definition \ref{D:operad}. 
The classical bar construction $\bar B^{GK}(\cal P)$ of an operad $\cal P$ is the cofree cooperad generated by $\Sigma\tilde{\cal P}$ with unique coderivation extending the partial composition  on $\cal P$ (see \cite[section 2]{GetJon94}). It has a description in terms of trees and it is graded by the number of vertices of the trees. 
Indeed, one has, for any finite set $I$ 

\begin{equation*}
 \bar B^{GK}_n(\cal P)(I)=K_{n-1}(b_I,\TI,\cal P).
 \end{equation*}

B. Fresse in \cite[section 4.1]{Fr04} builds also a complex $B^{GK}(R,\cal P,L)=R\circ \bar B^{GK}(\cal P)\circ L$  and proves that it satisfies the principal \ref{S:principle} of the  bar construction with coefficients.
 
He builds the levelization morphism

$$\Phi(R,\cal P,L):  B^{GK}(R,\cal P,L)\rightarrow N^{\circ}(R,\cal P,L)$$ and proves that it is a quasi-isomorphism.

In particular, for any finite set $I$, the quasi-isomorphism

\begin{align}
\Phi(\cal I,\cal P,\cal I)(I)_{n+1}: K_n(b_I,\TI,\cal P)=B^{GK}_{n+1}(\cal I,\cal P,\cal I)(I)\rightarrow N^{\circ}_{n+1}(\cal I,\cal P,\cal I)(I)
 \end{align}
 is described as follows.

 Let $t\in \TI$ be a tree with $n$ internal edges $:{e_1,\ldots, e_n}$. The source of an internal edge is the adjacent vertex closest to the leaves of the tree and its target is the adjacent vertex 
closest to the root of the tree. The set of internal edges of a tree $t$ is partially ordered: let $e$ and $f$ be internal edges, 
 $e\leq f$ if there is a path from a leaf of $t$ to the root of $t$ meeting $f$ before $e$. As an example the following  figure
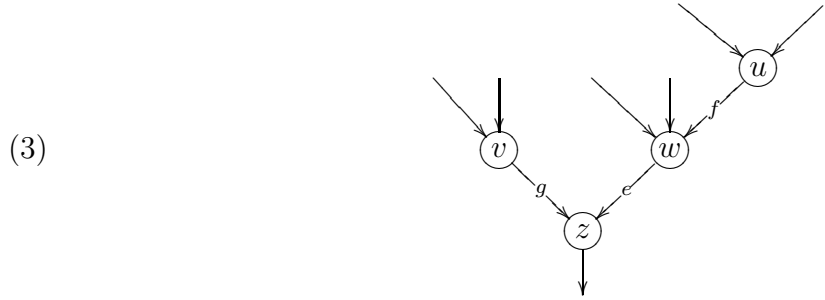
\begin{equation}\label{F:tree_efg}
\vcenter{\xymatrix@R=6mm@C=6mm{
	&&&\ar[dr]&&\ar[dl] \\
	\ar[dr]&\ar[d]&\ar[dr]&\ar[d]& *+<3mm>[o][F]{u}\ar[dl]|f& \\
	& *+<3mm>[o][F]{v}\ar[dr]|g &&  *+<3mm>[o][F]{w}\ar[dl]|e& \\
	&&  *+<3mm>[o][F]{z} \ar[d] &&& \\
	&&&&}}
\end{equation}
represents an element in $\cal P(t)$, where $\{e,f,g\}$ is the set of internal edges of $t$, with the partial order $e\leq f$. The source of $g$, $e$ and $f$ are $v$, $w$ and $u$ respectively. The target of $g$ and $e$ is $z$ and the target of $f$ is $w$.

Let  $e_1\wedge\ldots\wedge e_n\otimes p_t$ be an element in $\kfield (e_1\wedge\ldots\wedge e_n)\otimes \cal P(t)\subset K_n(b_I,\TI,\cal P)$. The levelization morphism associates to this element a sum of trees with $(n+1)$-levels.
The set $\{e_1,\ldots,e_n\}$ is partially ordered because it is the set of internal edges of a tree $t\in\TI$. The set $\{1,\ldots,n\}$ is totally ordered as a subset of $\N$.
To any  order-preserving bijection
$\sigma: \{e_1,\ldots,e_n\}\rightarrow \{1,\ldots,n\}$
one associates the level tree $t_\sigma$ where the source of $e_i$ is placed at level $\sigma(e_i)$, and where we complete the tree by adding vertices labelled by $1$ in $\cal P(1)$. The resulting element in $\cal P^{\circ (n+1)}$ is denoted by $\sigma(p_t)$.
The signature of $\sigma$, denoted by $\epsilon(\sigma)$ is the signature of the permutation $i\mapsto \sigma(e_i)$. The levelization morphism is defined by the following formula

\begin{equation}\label{E:levelization}
\Phi(\cal I,\cal P,\cal I)_{n+1}(I)(e_1\wedge\ldots\wedge e_n\otimes p_t)=
\sum_{\sigma: \{e_1,\ldots,e_n\}\rightarrow \{1,\ldots,n\} \text{ order-preserving }}
\epsilon(\sigma) \sigma(p_t)
\end{equation}

As an example we compute the levelization morphism associated to the element $e\wedge f\wedge g\otimes p_t$ of figure
(\ref{F:tree_efg}). The order-preserving maps involved in the formula (\ref{E:levelization}) are $\sigma_1:(e,f,g)\mapsto (1,2,3)$, $\sigma_2:(e,f,g)\mapsto (2,3,1)$ and
$\sigma_3:(e,f,g)\mapsto (1,3,2)$.

\begin{equation}
\begin{array}{ccccc}
\Phi(e\wedge f\wedge g\otimes p_t)= & & & & \\
& & & &  \\
\sigma_1(p_t)&+&\sigma_2(p_t)&-&\sigma_3(p_t)\\
& & & &  \\
\def\objectstyle{\scriptstyle}
\vcenter{\xymatrix@M=0pt@R=2mm@C=2mm{ 
 \ar[dr] && \ar[dl]  &\ar[d] & \ar[d] & \ar[d] &  \ar[d] \\
*+<3mm>{3}\ar@{.}[r] & *+<3mm>[o][F]{v}\ar@{.}[rr]\ar_g[d] & & *+<3mm>[o][F]{1}\ar@{.}[r]\ar[d] & 
*+<3mm>[o][F]{1}\ar@{.}[r]\ar[d] &  *+<3mm>[o][F]{1}\ar@{.}[r]\ar[d] &*+<3mm>[o][F]{1}\ar[dl] \\
*+<3mm>{2}\ar@{.}[r] & *+<3mm>[o][F]{1}\ar@{.}[rr]\ar_g[d] & & *+<3mm>[o][F]{1}\ar@{.}[r]\ar[dr] & 
 *+<3mm>[o][F]{1}\ar@{.}[r]\ar[d] & 
*+<3mm>[o][F]{u}\ar@{.}[r]\ar^{f}[dl]&  \\
*+<3mm>{1}\ar@{.}[rr] &  *+<3mm>[o][F]{1}\ar@{.}[rrr]\ar_g[dr] & & & *+<3mm>[o][F]{w}\ar@{.}[rr]\ar^e[dll] & &  \\
*+<3mm>{0}\ar@{.}[rr] &   &*+<3mm>[o][F]{z}\ar@{.}[rrrr]& &  & &  \\
 }}& +&  
 \def\objectstyle{\scriptstyle}\vcenter{\xymatrix@M=0pt@R=2mm@C=2mm{ & \ar[d] && \ar[d]  &\ar[d] & \ar[d] & \ar[d] &  \ar[dl] \\
*+<3mm>{3}\ar@{.}[r] & *+<3mm>[o][F]{1}\ar@{.}[rr]\ar[d] & & *+<3mm>[o][F]{1}\ar@{.}[r]\ar[d] & 
*+<3mm>[o][F]{1}\ar@{.}[r]\ar[dr] &*+<3mm>[o][F]{1}\ar@{.}[r]\ar[d]& *+<3mm>[o][F]{u}\ar@{.}[r]\ar^f[dl]& \\
*+<3mm>{2}\ar@{.}[r] & *+<3mm>[o][F]{1}\ar@{.}[rr]\ar[dr] & & *+<3mm>[o][F]{1}\ar@{.}[rr]\ar[dl] & 
&*+<3mm>[o][F]{w}\ar@{.}[rr]\ar^{e}[d]& & \\
*+<3mm>{1}\ar@{.}[rr] & &*+<3mm>[o][F]{v}\ar@{.}[rrr]\ar_g[dr] &  & 
&*+<3mm>[o][F]{1}\ar@{.}[rr]\ar^e[dll]& & \\
*+<3mm>{0}\ar@{.}[rrr] & & &*+<3mm>[o][F]{z}\ar@{.}[rrrr]& &  & &  \\
 }}  & - & \def\objectstyle{\scriptstyle}
\vcenter{\xymatrix@M=0pt@R=2mm@C=2mm{ 
& \ar[d] & \ar[d] &\ar[d] & \ar[d] & \ar[d] &  \ar[dl] \\
*+<3mm>{3}\ar@{.}[r] & *+<3mm>[o][F]{1}\ar@{.}[r]\ar[d]  & *+<3mm>[o][F]{1}\ar@{.}[r]\ar[dl] &  
*+<3mm>[o][F]{1}\ar@{.}[r]\ar[d] &*+<3mm>[o][F]{1}\ar@{.}[r]\ar[d]& *+<3mm>[o][F]{u}\ar@{.}[r]\ar^f[d]& \\
*+<3mm>{2}\ar@{.}[r] & *+<3mm>[o][F]{v}\ar@{.}[rr]\ar_g[d] & & *+<3mm>[o][F]{1}\ar@{.}[r]\ar[dr] 
&*+<3mm>[o][F]{1}\ar@{.}[r]\ar[d]& *+<3mm>[o][F]{1}\ar@{.}[r]\ar^{f}[dl]& \\
*+<3mm>{1}\ar@{.}[r]  &*+<3mm>[o][F]{1}\ar@{.}[rrr]\ar_g[dr]  &  & 
&*+<3mm>[o][F]{w}\ar@{.}[rr]\ar^e[dll]& & \\
*+<3mm>{0}\ar@{.}[rr]   & &*+<3mm>[o][F]{z}\ar@{.}[rrrr]& &  & &  \\
 }} 
\end{array}
\label{E:sigma_tree}
\end{equation}

\subsection{The factorization of $\overline\kappa:K(b_I,\TI,\cal P)\rightarrow N(b_I,\TI,\cal P)$}\label{S:main}\hfill\break
 
 We have seen in corollary \ref{C:kappa_qiso} that $\overline\kappa$
 is a quasi-isomorphism and that $B(b_I,\TI,\cal P)$ is identified with $B(\pi_I(b_I),\cal F,\cal P)$ in proposition \ref{P:evalfunctor}.

 The aim of this section is to prove that there exists a map (which will turn out to be a quasi-isomorphism)
 
 $$\overline\psi:N^{\circ}_{*}(\cal I,\cal P,\cal I)(I)\rightarrow N_{*-1}(b_I,\TI,\cal P)=N_{*-1}(\pi_I(b_I),\cal F,\cal P)$$
 
 such that $\overline\psi\Phi(\cal I,\cal P,\cal I)(I)=\overline\kappa$, that is, the following diagram is commutative
\begin{equation*}
\xymatrix{&K_{*-1}(b_I,\TI,\cal P)\ar[dl]_{\Phi(\cal I,\cal P,\cal I)(I)}\ar[dr]^{\overline\kappa}& \\
N^{\circ}_{*}(\cal I,\cal P,\cal I)(I)\ar[rr]_{\overline\psi}&&  N_{*-1}(b_I,\TI,\cal P) \\
    }
\end{equation*}

 We start with the description of a map
 $$\psi: B^{\circ}_{n+1}(\cal I,\cal P,\cal I)(I)=\underbrace{\cal P\circ\ldots\circ\cal P}_{n+1 \text{ terms }}(I)\rightarrow B_n(b_I,\TI,\cal P).$$
 
 An element $p_t$ in $B^{\circ}_{n+1}(\cal I,\cal P,\cal I)(I)$ is represented by a tree with $n+1$-levels, counted from $0$ to $n$ with vertices labelled by elements in $\cal P$. Such a level tree has subtrees of the form
 $$\xymatrix@M=0pt@R=2mm@C=2mm{ 
 *+<3mm>{l}\ar@{.}[r]  &  *+<3mm>[o][F]{w}\ar@{.}[r]\ar[d]& \\
  *+<3mm>{l-1}\ar@{.}[r]  &  *+<3mm>[o][F]{1}\ar@{.}[r]\ar[d]& \\
 *+<3mm>{\cdots}\ar@{.}[r]  & *+<3mm>[o][F]{1} \ar@{.}[r] \ar[d]& \\
  *+<3mm>{p+1}\ar@{.}[r]  &  *+<3mm>[o][F]{1}\ar@{.}[r]\ar[d]& \\
   *+<3mm>{p}\ar@{.}[r]  &  *+<3mm>[o][F]{z}\ar@{.}[r]& \\
    }$$
with $w\in\cal P(x), x\geq 2$ and $z\in\cal P(y),y\geq 2$. We define the {\it level-edge set} $N(p_t)$ as $e\in N(p_t)$ if and only if 
there is a sequence of consecutive edges in
$p_t$, $e=\{e_1>\ldots>e_k\}$ such that the source of $e_1$ lives in $\cal P(x)$ with $x\geq 2$, the target of $e_k$ lives in $\cal P(y)$ with $y\geq 2$ and all other sources and targets lives in $\cal P(1)$. The source of $e$ is the source of $e_1$ and the target of $e$
is the target of $e_k$. The levels of the source and target of $e$ are denoted by  $s(e)$ and $t(e)$ respectively.  The previous figure shows an element $e$ in $N(p_t)$ such that $s(e)=l$ and $t(e)=p$. The idea underlying the definition of $N(p_t)$ is that we don't want to consider vertices labelled by $1\in\cal P(1)$.   
For $1\leq i\leq n$, let
 $$N_i(p_t)=\{e\in N(p_t)| t(e)<i\leq s(e)\}.$$
 One has
 $N(p_t)=\cup_{1\leq i\leq n} N_i(p_t),$
 for $0\leq t(e)< n$ and $1\leq s(e)\leq n$. Note that this decomposition is not necessarily a partition of $N(p_t)$ as we will see in example \ref{Ex:N}.
 Let $t$ be a level tree and $p_t\in\otimes_{v\in V_t} \cal P(\In(v))$. By forgetting the units, we denote by $r(t)$ the associated rooted tree and by $r(p_t)$ the associated element in $\cal P(r(t))$.  In the sequel the level-edge set $N(p_t)$ is written according to its decomposition $N(p_t)=(N_1(p_t),\ldots,N_n(p_t))$.

\begin{sub}{\it Example:}\label{Ex:N}  The associated reduced tree to any of the trees of equation (\ref{E:sigma_tree}) is the tree
$$\vcenter
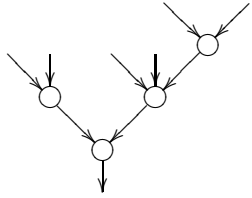
{\xymatrix@R=4mm@C=4mm{
	&&&\ar[dr]&&\ar[dl] \\
	\ar[dr]&\ar[d]&\ar[dr]&\ar[d]& *+<3mm>[o][F]{}\ar[dl]& \\
	& *+<3mm>[o][F]{}\ar[dr]&&  *+<3mm>[o][F]{}\ar[dl]& \\
	&&  *+<3mm>[o][F]{} \ar[d] &&& \\
	&&&&}}$$
 and the associated element $r(p_t)$ is the tree of figure  (\ref{F:tree_efg}). The set of level-edges of
$$q_t=\def\objectstyle{\scriptstyle}
\vcenter
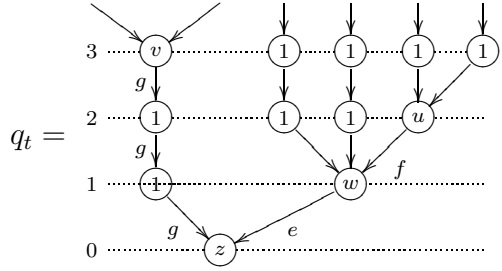
{\xymatrix@M=0pt@R=4mm@C=4mm{ 
 \ar[dr] && \ar[dl]  &\ar[d] & \ar[d] & \ar[d] &  \ar[d] \\
*+<3mm>{3}\ar@{.}[r] & *+<3mm>[o][F]{v}\ar@{.}[rr]\ar_g[d] & & *+<3mm>[o][F]{1}\ar@{.}[r]\ar[d] & 
*+<3mm>[o][F]{1}\ar@{.}[r]\ar[d] &  *+<3mm>[o][F]{1}\ar@{.}[r]\ar[d] &*+<3mm>[o][F]{1}\ar[dl] \\
*+<3mm>{2}\ar@{.}[r] & *+<3mm>[o][F]{1}\ar@{.}[rr]\ar_g[d] & & *+<3mm>[o][F]{1}\ar@{.}[r]\ar[dr] & 
 *+<3mm>[o][F]{1}\ar@{.}[r]\ar[d] & 
*+<3mm>[o][F]{u}\ar@{.}[r]\ar^{f}[dl]&  \\
*+<3mm>{1}\ar@{.}[rr] &  *+<3mm>[o][F]{1}\ar@{.}[rrr]\ar_g[dr] & & & *+<3mm>[o][F]{w}\ar@{.}[rr]\ar^e[dll] & &  \\
*+<3mm>{0}\ar@{.}[rr] &   &*+<3mm>[o][F]{z}\ar@{.}[rrrr]& &  & &  \\
 }}$$
is $N(q_t)=\{e,f,g\}$ with $N_1(q_t)=\{e,g\}, N_2(q_t)=\{f,g\}$ and 
$N_3(q_t)=\{g\}$.
\end{sub}

\medskip
  
\begin{defn} Let $p_t$ be a level tree and $N$ its associated level-edge set.
For $\sigma\in S_n$, we define  $E_1^\sigma=N_{\sigma(1)}$ and $E_i^\sigma=N_{\sigma(i)}\setminus \{ \petitcup\limits_{1\leq j\leq i-1} N_{\sigma(j)}\}=N_{\sigma(i)}\setminus \{ \petitcup\limits_{1\leq j\leq i-1} E_j^\sigma\}$. The map $\psi$ is defined as follows

\begin{equation*}
 \begin{array}{cccc}
\psi: &  B^{\circ}_{n+1}(\cal I,\cal P,\cal I)(I)&\rightarrow &B_n(b_I,\TI,\cal P) \\
&p_t & \mapsto & \sum_{\sigma\in S_n}\epsilon(\sigma) (E_1^\sigma,\ldots,E_n^\sigma)\otimes r(p_t),\\
\end{array} 
\end{equation*}

\end{defn}

\medskip
\noindent{\it Example.} The computation of $\psi(q_t)$ for the tree $q_t$ of example  \ref{Ex:N} gives
$$\psi(q_t)=(\underbrace{(\{e,g\},f,\emptyset)}_{\sigma=(123)}-
\underbrace{(\{f,g\},e,\emptyset)}_{\sigma=(213)}+\underbrace{(\{f,g\},\emptyset,e)}_{\sigma=(231)}-
\underbrace{(g,f,e)}_{\sigma=(321)}+\underbrace{(g,e,f)}_{\sigma=(312)}-\underbrace{(\{e,g\},\emptyset,f)}_{\sigma=(132)})\otimes r(p_t).$$

\begin{lem} The map $\psi$ induces a well-defined map
$$\overline\psi:   N^{\circ}_{n+1}(\cal I,\cal P,\cal I)(I)\rightarrow N_n(b_I,\TI,\cal P),$$
which commutes with the differentials.
 \end{lem}
 
\begin{proof} Assume $p_t=s_j(q)$ with $s_j:  B^{\circ}_{n}(\cal I,\cal P,\cal I)(I)\rightarrow  B^{\circ}_{n+1}(\cal I,\cal P,\cal I)(I)$ being the degeneracy map sending $\cal I\circ \cal P^{\circ n}\circ \cal I$ to $\cal I\circ\cal P^{\circ j}\circ \cal I\circ \cal P^{\circ n-j}\circ \cal I$. If 
$1\leq j\leq n-1$, then the vertices of the tree $t$ at level $j$ are all labelled by $1$.
Consequently, $N_j(p_t)=N_{j+1}(p_t)$ and using the transposition $(j\, j+1)$ one gets that $\psi(p_t)=0$. If $j=0$, then $N_1(p_t)=\emptyset$ and the composite of $\psi$ with the projection $B_n(b_I,\TI,\cal P)\rightarrow N_n(b_I,\TI,\cal P)$ is zero. If $j=n$ then $N_n(p_t)=\emptyset$ and the composite of $\psi$ with the projection $B_n(b_I,\TI,\cal P)\rightarrow N_n(b_I,\TI,\cal P)$ is zero.

\smallskip 
 
In order to prove that for every $x$ in  $N^{\circ}_{n+1}(\cal I,\cal P,\cal I)(I)$ one has $\overline\psi(dx)=d\overline\psi(x)$, it is enough to prove the equality for a representant $p_t$ of $x$ in $B^{\circ}_{n}(\cal I,\cal P,\cal I)(I)$,  
such that $N_j(p_t)\not=\emptyset,\forall j$. To keep track of the levels we write such an
element $(N_1,\ldots,N_n,p_t)$, where we consider  $p_t\in\cal P(r(t))$, forgetting the units. 
On the one hand
the differential is given by
$$d(N_1,\ldots,N_n,p_t)=\sum_{i=1}^n(-1)^i \left(N_1,\ldots,\hat N_i,\ldots, N_n,(N_i\setminus{\cup_{j\not=i} N_j})_*(r(p_t))\right).$$
Identifying permutations in $S_{n-1}$ with permutations $\sigma$ in $S_n$ such that $\sigma(n)=i$ one gets
$$\overline\psi d(N_1,\ldots,N_n,r(p_t))=\sum_{i=1}^n(-1)^n\sum_{\sigma\in S_n| \sigma(n)=i}\epsilon(\sigma) (E_1^\sigma,\ldots,E_{n-1}^\sigma)\otimes 
(N_i\setminus{\cup_{j\not=i} N_j})_*(p_t).$$
On the other hand one has
\begin{multline*}
d\overline\psi (N_1,\ldots,N_n,p_t)=\sum_{i=1}^n(-1)^i d_i\sum_{\sigma\in S_n}
\epsilon(\sigma)(E_1^\sigma,\ldots,E_n^\sigma)\otimes p_t\\
=(-1)^n\sum_{\sigma\in S_n} \epsilon(\sigma)(E_1^\sigma,\ldots,E_{n-1}^\sigma)\otimes (E_n^\sigma)_*(p_t),
\end{multline*}
 for regrouping permutations by pairs $(\sigma,\tau)$ such that $\sigma(i)=k,\sigma(i+1)=l$ and
$\tau(i)=l,\tau(i+1)=k$, one gets  $d_i(\sum_{\sigma\in S_n}
(E_1^\sigma,\ldots,E_n^\sigma)\otimes r(p_t))=0$ if $1\leq i<n$.

Furthermore $E_n^{\sigma}=N_{\sigma(n)}\setminus{\cup_{j\not=\sigma(n)}N_j}$ implies 
\begin{multline*}
d\overline\psi (N_1,\ldots,N_n,p_t)=
\sum_{i=1}^n(-1)^n\sum_{\sigma\in S_n| \sigma(n)=i}\epsilon(\sigma) (E_1^\sigma,\ldots,E_{n-1}^\sigma)\otimes 
(N_i\setminus{\cup_{j\not=i} N_j})_*(p_t).
\end{multline*}
 The two expressions coincide.
\end{proof}

\begin{thm}The map $\overline\kappa:K(b_I,\TI,\cal P)\rightarrow N(b_I,\TI,\cal P)$ factorizes through  $N_*^\circ(\cal I,\cal P,\cal I)(I)$, and the following diagram
\begin{equation*}
\xymatrix{&K_{*-1}(b_I,\TI,\cal P)\ar[dl]_{\Phi(\cal I,\cal P,\cal I)(I)}\ar[dr]^{\overline\kappa}& \\
N^{\circ}_{*}(\cal I,\cal P,\cal I)(I)\ar[rr]_{\overline\psi}&&  N_{*-1}(b_I,\TI,\cal P) \\
    }
\end{equation*}
is commutative. Consequently  $\overline\psi$ is a quasi-isomorphism.
 \end{thm}

\begin{proof} The symbol $[k]$ denotes the set $\{1,\ldots,k\}$. Recall that
$$\overline\psi\Phi(e_1\wedge\ldots\wedge e_n\otimes p_t)=\overline\psi(\sum_{f:\{e_1,\ldots,e_n\}\rightarrow \{1,\ldots,n\}\atop{ \text{ order-preserving}} }
\epsilon(f) f(p_t)).$$
We prove the theorem by induction on $n$. If $n=1$ it is obvious. Assume the result is true for any tree with $n-1$ internal edges.
Let $p_t$ be a tree with $n$ internal edges $E=\{e_1,\ldots,e_n\}$. One can re-order the internal edges so that, there is a chain of consecutive edges from the root to a leaf $a_1<\ldots <a_p$ such that $a_p=e_n$ and if $p>1$, then $a_{p-1}=e_{n-1}$.  By convention, if $p=1$, we let $a_{0}=\emptyset$. 

Let $\tilde p_t$ be the tree obtained from $p_t$ by removing the edge $e_n$. It has exactly $n-1$ internal edges $\tilde E=\{e_1,\ldots,e_{n-1}\}$. Let 
$f:\tilde E\rightarrow [n-1]$ be an order-preserving map. By convention $f(\emptyset)=0$.
For  $f(a_{p-1})<i\leq n$, let us define
$$\begin{array}{lccc}
f^i:&\{e_1,\ldots,e_n\}&\rightarrow& [n] \\
& e_j, \ j<n & \mapsto& \begin{cases}  
           f(e_j) & \text{ if } f(e_j)<i \\
            f(e_j)+1 & \text{ if } f(e_j)\geq i \\
          \end{cases}\\
& e_n &\mapsto & i\end{array}$$
The map $f^i$ is an order-preserving bijection. One has $\epsilon(f^n)=\epsilon(f)$. As a consequence $\epsilon(f^i)=(-1) ^{n-i}\epsilon(f)$, for  $f^i=(i\, \ldots \, n)f^n$ where $(i\, \ldots \, n)$ denotes the cycle
$i\rightarrow i+1\rightarrow\ldots\rightarrow n\rightarrow i$ in $S_n$.

Furthermore, to any order-preserving bijection
$\tau:E\rightarrow [n]$ there exists a unique $f:\tilde E\rightarrow [n-1]$ and a unique
$i$ with $f(a_{p-1})<i\leq n$  such that $f^i=\tau$. 

Consequently

$$\Phi(e_1\wedge\ldots\wedge e_n\otimes p_t)=\sum_{f:\tilde E\rightarrow [n-1]\atop{ \text{ order-preserving}}} \sum_{i=f(a_{p-1})+1}^n
(-1)^{n-i}\epsilon(f) f^i(p_t).$$

In order to evaluate $\overline\psi$ on the above expression, one needs to express $N_k^{f^i}(p_t)$ in terms of the $N_j^{f}(\tilde p_t)$'s. Because of the choices of the level for $f^i(p_t)$, one has, for $f(a_{p-1})<i\leq n$,

\begin{equation}\label{E:N_sigma_i}
N_k^{f^i}(p_t)=\begin{cases} N_k^f(\tilde p_t), & \text{ if } k\leq f(a_{p-1})<i, \\
N_k^f(\tilde p_t)\cup e_n, & \text{ if } f(a_{p-1}) <k\leq i, \\
N_{k-1}^f(\tilde p_t),& \text{ if } i<k\leq n.
\end{cases}
\end{equation}

Note that if $i=n$ the second equality reads $N_n^{f^n}(p_t)=\{e_n\}$ since $N_n^f=\emptyset$. 

For example, if $n=4$ and $f(a_{p-1})=1$, writing the sets $N^{f^i}$ as
$(N^{f^i}_1,N^{f^i}_2,N^{f^i}_3,N^{f^i}_4)$ one gets
\begin{align*}
N^{f^2}=&(N^f_1,N^f_2\cup e_4,N^f_2,N^f_3) \\
N^{f^3}=&(N^f_1,N^f_2\cup e_4,N^f_3\cup e_4,N^f_3) \\
N^{f^4}=&(N^f_1,N^f_2\cup e_4,N^f_3\cup e_4,e_4). \\
\end{align*}

Recall that 
$$\bar\psi(f(p_t))=\sum_{\sigma\in S_n} \epsilon(\sigma)(E^{f,\sigma}_1,\ldots,E^{f,\sigma}_n)\otimes p_t,
\text{ with } E^{f,\sigma}_k=N^f_{\sigma(k)}\setminus {\petitcup\limits_{i<k} N^f_{\sigma(i)}}.$$
Let $\sigma\in S_n$, $j=f(a_{p-1})$ and $j<i\leq n$. By relations (\ref{E:N_sigma_i}), the set $N^{f^i}$  decomposes as
$$N^{f^i}=(N_1^f,\ldots,N^f_j,N^f_{j+1}\cup \{e_n\},\ldots,N^f_i\cup \{e_n\},N^f_i,N^f_{i+1},\ldots,N_{n-1}^f).$$

Firstly, if $\sigma^{-1}(i)=k<\sigma^{-1}(i+1)=l$ then the sequence
$(E^{f^i,\sigma}_1,\ldots,E^{f^i,\sigma}_n)$ satisfies $E^{f^i,\sigma}_l=\emptyset$ and vanishes in $N_n(b_I,\TI,\cal P)$. Hence we only need to consider the elements $\sigma\in S_n$ such that $\sigma^{-1}(i+1)<\sigma^{-1}(i)$. In that case $E^{f^i,\sigma}_{\sigma^{-1}(i)}=\{e_n\}$.

Secondly, for $j+1\leq r\leq i-1$, if  $\sigma^{-1}(r)=k<\sigma^{-1}(i)=l$ then the sequence
$(E^{f^i,\sigma}_1,\ldots,E^{f^i,\sigma}_n)$ satisfies $E^{f^i,\sigma}_l=\emptyset$ and vanishes in $N_n(b_I,\TI,\cal P)$. 

As a consequence, we only need to consider  the elements $\sigma\in S_n$ such that 
$$\sigma^{-1}(i+1)<\sigma^{-1}(i)<\{\sigma^{-1}(j+1),\ldots,\sigma^{-1}(i-1)\}.$$

Note that if $i=n$ the latter condition writes
$$\sigma^{-1}(n)<\{\sigma^{-1}(j+1),\ldots,\sigma^{-1}(n-1)\}.$$
and if $i=j+1$ it writes
$$\sigma^{-1}(i+1)<\sigma^{-1}(i).$$

Let $1\leq l\leq n$ be a fixed integer. Choose $\sigma\in S_n$ such that $\sigma^{-1}(i)=l$. The condition  $\sigma^{-1}(i+1)<\sigma^{-1}(i)<\{\sigma^{-1}(j+1),\ldots,\sigma^{-1}(i-1)\}$ implies that the sequence
$(E^{f^i,\sigma}_1,\ldots,E^{f^i,\sigma}_n)$ writes $(E^{f,\tau}_1,\ldots,E^{f,\tau}_{l-1},\{e_n\},E^{f,\tau}_{l+1},\ldots,E^{f,\tau}_{n-1})$ with $\tau\in S_{n-1}$ obtained as the composite
$\sigma_i\sigma\delta_l$ where $\delta_l:[n-1]\rightarrow [n]$ is the map missing $l$ and $\sigma_i :[n]\rightarrow [n-1]$ is the map repeating $i$. It is clear that $\epsilon(\sigma)=\epsilon(\tau)(-1) ^{l+i}$.
When $i$ runs from $j+1$ to $n$ one covers $S_{n-1}$. For, if $i=j+1$ then the set involving $N_{j+1}^f$ appears before $e_n$ and if $i>j+1$ then it appears after $e_n$. If $i=j+2$ then
 the set involving $N_{j+2}^f$ appears before $e_n$ and if $i>j+2$ then it appears after $e_n$. And so on.

It yields the computation:
\begin{multline*}
\bar\psi\Phi(e_1\wedge\ldots\wedge e_n\otimes p_t)=\sum_{f:\tilde E\rightarrow [n-1]\atop{ \text{ order-preserving}}} \sum_{i=f(a_{p-1})+1}^n
(-1)^{n-i}\epsilon(f)\sum_{l=1}^{n} \sum_{\sigma\in S_n,\atop{ \sigma^{-1}(i)=l}} \epsilon(\sigma)(E^{f,\sigma}_1,\ldots,E^{f,\sigma}_n)\otimes p_t= \\
\sum_{f:\tilde E\rightarrow [n-1]\atop{ \text{ order-preserving}}}
\epsilon(f)\sum_{l=1}^n \sum_{\sigma\in S_{n-1}} (-1)^{n+l}\epsilon(\sigma)(E^{f,\sigma}_1,\ldots,E^{f,\sigma}_{l-1},e_n,E^{f,\sigma}_l,\ldots,E^{f,\sigma}_{n-1})\otimes p_t=\\
\sum_{l=1}^n
\sum_{\sigma\in S_{n-1}} (-1)^{n+l}\epsilon(\sigma)(e_{\sigma(1)},\ldots,e_{\sigma(l-1)},e_n,e_{\sigma(l+1)},\ldots,e_{\sigma(n-1)})\otimes p_t=\\
\sum_{\sigma\in S_{n}} \epsilon(\sigma)(e_{\sigma(1)},\ldots,e_{\sigma(n)})\otimes p_t=
\bar\kappa(e_1\wedge\ldots\wedge e_n\otimes p_t).
\end{multline*}
\end{proof}

 \section{Resolution of the category $\TI$ and operads up to homotopy}

This section is devoted to the bar and cobar construction for differential graded categories and cocategories whose objects are the objects of $\TI$. It follows closely the paper  \cite{Frpr09}. In this paper B. Fresse works on the category of Batanin trees $\Epi_n$. He proves that the complex we obtained with B. Richter in \cite{LivRic09pr} corresponds to the bar construction of the category $\Epi_n$ with coefficients in the Loday functor and the unit functor. In his paper, he proves that the category $\Epi_n$ is Koszul, yielding a minimal model $R(\Epi_n)$ of $\Epi_n$.

In this paper we work exactly in the same spirit; we have proved in section \ref{SS:cat} that the category $\TI$ is Koszul and the purpose
of the first section is to express it's minimal model $R(\TI)\rightarrow \TI$. The main result will be that given a species $\cal M$, then the map $R(\TI)\rightarrow \dgvs$ which associates $\cal M(t)$ to $t\in\TI$ is a functor if and only if $\cal M$ is an operad up to homotopy.

\subsection{Bar and cobar construction for dg categories and dg cocategories}\hfill\break

The bar and cobar constructions follows closely the ones for associative and coassociative algebras, and in this section we just state our notation and the theorem needed for the sequel.

\medskip

From now on we denote by $\Ob\TI$ the set of trees in the category $\TI$. A tree $t$ has a degree $|t|$ given by the number of internal edges. A {\sl dg graph} is a map $\Gamma:\Ob\TI\times \Ob\TI\rightarrow \dgvs$.
We denote by  $\cal C_I$ the category of differential graded connected categories whose objects are the trees 
 $t\in \Ob\TI$. Let $\cal C$ be such a category. Such a data is equivalent to
\begin{itemize}
 \item A dg graph $\cal C$ which will corresponds to the morphisms in the category.
 \item For every $a,b,c\in\Ob\TI$, composition maps $\cal C(b,c)\otimes \cal C(a,b)\rightarrow \cal C(a,c)$ in $\dgvs$ which are associative; 
\item Identity elements  $1_a\in \cal C(a,a)$ which are unit for the composition;
 \item Connectivity assumption: $\forall a\in \Ob\TI, \cal C(a,a)=\kfield$ and
$\cal C(b,a)=0$ if $|b|<|a|$.
\end{itemize}
An example of such a category is $\kfield\TI$.

For a connected dg graph $\Gamma$,
we denote by $\overline{\Gamma}(s,t)=\begin{cases} 0, & \text{ if } s=t, \\ \Gamma(s,t), & \text{ if } s\not=t.\end{cases}$

Similarly we define $\cal {C}^c_I$ the category of differential graded connected cocategories whose objects are $\Ob\TI$. A cocategory is defined the same way as a category except that the arrows go in the reverse order.

Given a dg graph $\Gamma:\Ob\TI\times \Ob\TI\rightarrow \dgvs$, one can form the free category generated by $\Gamma$. As a dg graph, one has 

$$\cal F(\Gamma)(t,t')=\bigoplus\limits_{t'=x_0,\ldots,x_m=t} \Gamma(x_1,x_0)\otimes \Gamma(x_2,x_1)\otimes\ldots\Gamma(x_{m-1},x_{m-2})\otimes\ldots\Gamma(x_m,x_{m-1}).$$
The compositions of maps are given by the concatenation.

Similarly the free co-category generated by $\Gamma$, denoted by $\cal F^c(\Gamma)$ is given by the same dg graph and the co-composition  are given by the deconcatenation.

There is an adjunction between co-categories and categories

$$\Omega: \text{cocategories} \rightleftarrows \text{categories} : B$$

The bar construction $B(\cal C)$ of the category $\cal C$ is the free cocategory $\cal F^c(s\overline{\cal C})$ with the unique coderivation lifting the composition product in $\cal C$. Namely   
$$\partial(s\alpha_1\otimes\ldots\otimes s\alpha_p)=
\sum_{i=1}^{p-1}(-)^{|s\alpha_1|+\ldots |s\alpha_i|} s\alpha_1\otimes\ldots\otimes
s(\alpha_i\alpha_{i+1})\otimes\ldots s\alpha_p,$$
where $|s\alpha_i|$ denotes the degree of $s\alpha_i\in s\cal C(a,b)$ (see e.g. \cite{EML53})

The cobar construction $\Omega(\cal R)$ of a cocategory is the free category $\cal F(s^{-1}\overline{\cal R})$ with the unique derivation lifting the co-composition product in $\cal R$:  
$$\partial(s^{-1}\alpha_1\otimes\ldots\otimes s^{-1}\alpha_p)=\sum_{i=1}^{p}(-)^{|\alpha_1|+\ldots |\alpha_i|+i} 
s^{-1}\alpha_1\otimes\ldots\otimes 
s^{-1}\alpha_{i,(1)}\otimes s^{-1}\alpha_{i,(2)}\otimes\ldots s^{-1}\alpha_p,$$
where we use the Sweedler's notation for the co-composition.

\begin{lem} Let $\cal C$ be a category in $\cal C_I$. The counit of the adjunction is a quasi-isomorphism:
$$\Omega B(\cal C)\rightarrow \cal C.$$
\end{lem}

\subsection{Minimal model of the category $\TI$}

In this section,
we use the Koszul complex in order to build the minimal resolution of $\TI$.

\begin{lem}The bar construction $B(\kfield\TI)(t,s)$ corresponds to the normalized bar construction
$N(b_s,\TI,b_t)$.
\end{lem}

\begin{proof} The two definitions coincide, and we just have to check that the degrees and differentials coincide. For $t,s\in\Ob\TI$ with 
$s=t/E$, an element in 
$N_n(b_s,\TI,b_t)$ writes $(E_1,\ldots,E_n)$, with $E_1\sqcup E_2\sqcup\ldots\sqcup E_n=E$ and $E_i$ is non empty for every $i$. It has degree $n$ and its differential is given by
$d(E_1,\ldots,E_n)=\sum_{i=1}^{n-1} (-1)^i (E_1,\ldots,E_i\sqcup E_{i+1},\ldots, E_n)$.
\end{proof}

Consequently the dg graph $(t,s)\mapsto N(b_s,\TI,b_t)$ is endowed with a structure of cocategory. Note that for 
$s=t$ one has $N(b_t,\TI,b_t)$ is 1-dimensional concentrated in degree $0$.

Recall that 
$$K_*(b_s,\TI,b_t)=\begin{cases}\Lambda^{|E|}(\kfield[E]), & \text{ if } s=t/E  \text{ and } *=|E|,\\
0,&\text{ elsewhere}, \end{cases}$$
with zero differential. Hence $K$ determines a dg graph
$$\begin{array}{cccc}
K: & \Ob\TI\times\Ob\TI & \rightarrow& \dgvs \\
& (t,s) &\mapsto & K(b_s,\TI,b_t)
\end{array}$$

 We define the co-composition on $K$ by
$$\Delta(h=e_1\wedge\ldots\wedge e_n)=1\otimes h+h\otimes 1+\sum_{p=1}^{n-1} 
\sum_{\sigma\in\Sh_{p,n-p}} \epsilon(\sigma) e_{\sigma(1)}\wedge\ldots\wedge e_{\sigma(p)}\otimes  e_{\sigma(p+1)}\wedge\ldots\wedge e_{\sigma(n)},$$
where $\Sh_{p,n-p}$ denotes the set of $(p,n-p)$-shuffles.

\begin{lem}The   dg graph $K$ is  a subcocategory of $B(\kfield\TI)$ via the map $\bar\kappa:
K(b_s,\TI,b_t)\rightarrow N(b_s,\TI,b_t).$
\end{lem}

\begin{proof}The fact that the co-composition commutes with $\bar\kappa$ comes from the bijection between
$(S_p\times S_{n-p})\Sh_{p,n-p}$ and $S_n$.
\end{proof}

Because the category $\TI$ is Koszul (see theorem \ref{T:TIKoszul}) the morphism of co-categories $\bar\kappa$ is a quasi-isomorphism. Since $\Omega$ behaves well with respect to these quasi-isomorphisms, one has

\begin{thm} The cobar construction of the co-category $K$ is a minimal resolution of the category
$\kfield\TI$. 
\end{thm}

\subsection{Operads up to homotopy}

Let $\cal M$ be a vector species, and consider the map

$$\begin{array}{cccc}
\underline{\cal M}: &\Ob\TI & \rightarrow & \dgvs \\
& t & \mapsto & \cal M(t)=\bigotimes\limits_{v\in E_t} \cal M(\In(v)),\\
\end{array}$$
defined in section \ref{S:left_TI}.

\begin{thm} Let $\cal M$ be a vector species. The map $\underline{\cal M}$ determines a functor
$\Omega(K)\rightarrow \dgvs$ if and only if $\cal M$ is an operad up to homotopy.
\end{thm}

\begin{proof} Recall that $\Omega(K)=\cal F(s^{-1}K)$.

Assume that $\underline{\cal M}$ is a functor. Since $\Omega K$ is the free category generated by $s^{-1}K$, one has for every $t,s=t/E$ a composition map
$$\circ_E:\cal M(t)\rightarrow \cal M(s),$$
of degree $|E|-1$. 
Let us write $E=e_1\wedge\ldots\wedge e_n$  a generator of the one dimensional vector space $K(b_t,\TI,b_I)$. In $\Omega(K)$ one has
$$d(s^{-1}E)=\sum_{p=1}^{n-1}
\sum_{\sigma\in\Sh_{p,n-p}} \epsilon(\sigma)
s^{-1}(e_{\sigma(1)}\wedge\ldots\wedge e_{\sigma(p)})\otimes 
s^{-1}(e_{\sigma(p+1)}\wedge\ldots\wedge e_{\sigma(n)}).$$
In terms of functors, it writes
$$\partial(\circ_E)=\sum_{F\sqcup G=E, \atop{F,G\not=\emptyset}}\epsilon(F,G) \circ_F\circ_G,$$
where  $\epsilon(F,G)=\epsilon(\sigma)$ for  the shuffle $\sigma$ corresponding to the sets $F$ and $G$ when an order of elements in $E$ is given.
This is exactly the definition of an operad up to homotopy in \cite[4.2.2]{vanderlaanphD}.
\end{proof}

\bibliographystyle{amsplain}
\bibliography{bibliojanv2011}
\bigskip

\end{document}